\def\Version{{6}}




\message{<< Assuming 8.5" x 11" paper >>}    

\magnification=\magstep1	          
\raggedbottom

\parskip=9pt

%

\def\singlespace{\baselineskip=12pt}      
\def\sesquispace{\baselineskip=16pt}      

 \let\miguu=\footnote
 \def\footnote#1#2{{$\,$\parindent=9pt\baselineskip=13pt%
 \miguu{#1}{#2\vskip -7truept}}}
%
%

\def\linebreak{\hfil\break}
\def\lbr{\linebreak}
\def\pagebreak{\vfil\break}


\def\BulletItem #1 {\item{$\bullet$}{#1 }}
\def\bulletitem #1 {\BulletItem{#1}}

\def\PrintVersionNumber{
 \vskip -1 true in \medskip 
 \rightline{version \Version} 
 \vskip 0.3 true in \bigskip \bigskip}

\def\author#1 {\medskip\centerline{\it #1}\bigskip}

\def\address#1{\centerline{\it #1}\smallskip}

\def\furtheraddress#1{\centerline{\it and}\smallskip\centerline{\it #1}\smallskip}

\def\email#1{\smallskip\centerline{\it address for email: #1}} 

\def\AbstractBegins
{
 \singlespace                                        
 \bigskip\leftskip=1.5truecm\rightskip=1.5truecm     
 \centerline{\bf Abstract}
 \smallskip
 \noindent	
 } 
\def\AbstractEnds
{
 \bigskip\leftskip=0truecm\rightskip=0truecm       
 }

\def\section #1 {\bigskip\noindent{\headingfont #1 }\par\nobreak\noindent}

\def\subsection #1 {\medskip\noindent{\subheadfont #1 }\par\nobreak\noindent}
 %

\def\ReferencesBegin
{
 \singlespace					   
 \vskip 0.5truein
 \centerline           {\bf References}
 \par\nobreak
 \medskip
 \noindent
 \parindent=2pt
 \parskip=6pt			
 }
 %


\def\reference{\hangindent=1pc\hangafter=1} 

\def\ref{\reference}

 %

\def\journaldata#1#2#3#4{{\it #1\/}\phantom{--}{\bf #2$\,$:} $\!$#3 (#4)}
 %

\def\eprint#1{{\tt #1}}

\def\arxiv#1{\hbox{\tt http://arXiv.org/abs/#1}}
 %

\def\webtilde{\lower2pt\hbox{${\widetilde{\phantom{m}}}$}}

\def\webhome{{\tt {http://www.physics.syr.edu/}{\webtilde}{sorkin/}}}
 %

 %


\def\hpf#1{\webhome{\tt{some.papers/}}}
 %

\def\hpfll#1{\webhome{\tt{lisp.library/}}}
 %




\font\titlefont=cmb10 scaled\magstep2 

\font\headingfont=cmb10 at 12pt
%

\font\subheadfont=cmssi10 scaled\magstep1 
%




\def\QuotationBegins
{
 \singlespace
 \smallskip\leftskip=1.5truecm\rightskip=1.5truecm  
 \noindent	    
 } 
\def\QuotationEnds
{
 \smallskip\leftskip=0truecm\rightskip=0truecm      
 \sesquispace			
 \par				
 \noindent
 } 





\def\catuskoti{catu{\d s}ko{\d t}i }
\def\Catuskoti{Catu{\d s}ko{\d t}i }
\def\paryudasa{paryud{\= a}sa } 
\def\Nagarjuna{N{\= a}g{\= a}rjuna }

\def\lbrace{$\{$}
\def\rbrace{$\}$}

\def\QuotationBegins
{
 \singlespace                                
 \smallskip\leftskip=1.5truecm\rightskip=1.5truecm     
 \noindent	
 } 
\def\QuotationEnds
{
 \smallskip\leftskip=0truecm\rightskip=0truecm       
 \sesquispace
 \par				
 \noindent
 } 

\newcount\tmpnum \newdimen\tmpdim
{\lccode`\?=`\p \lccode`\!=`\t  \lowercase{\gdef\ignorept#1?!{#1}}}

\edef\widecharS{\expandafter\ignorept\the\fontdimen1\textfont1}

\def\widebar#1{\futurelet\next\widebarA#1\widebarA}
\def\widebarA#1\widebarA{%
   \def\tmp{0}\ifcat\noexpand\next A\def\tmp{1}\fi
   \widebarE
   \ifdim\tmp pt=0pt \overline{#1}%
   \else {\mathpalette\widebarB{#1}}\fi
}
\def\widebarB#1#2{%
   \setbox0=\hbox{$#1\overline{#2}$}%
   \tmpdim=\tmp\ht0 \advance\tmpdim by-.4pt
   \tmpdim=\widecharS\tmpdim
   \kern\tmpdim\overline{\kern-\tmpdim#2}%
}

\def\widebarC#1#2 {\ifx#1\end \else 
   \ifx#1\next\def\tmp{#2}\widebarD 
   \else\expandafter\expandafter\expandafter\widebarC
   \fi\fi
}
\def\widebarD#1\end. {\fi\fi}

\def\widebarE{\widebarC A1.4 J1.2 L.6 O.8 T.5 U.7 V.3 W.1 Y.2 
   a.5 b.2 d1.1 h.5 i.5 k.5 l.3 m.4 n.4 o.6 p.4 r.5 t.4 v.7 w.7 x.8 y.8
   \alpha1 \beta1 \gamma.6 \delta.8 \epsilon.8 \varepsilon.8 \zeta.6 \eta.4
   \theta.8 \vartheta.8 \iota.5 \kappa.8 \lambda.5 \mu1 \nu.5 \xi.7 \pi.6
   \varpi.9 \rho1 \varrho1 \sigma.7 \varsigma.7 \tau.6 \upsilon.7 \phi1
   \varphi.6 \chi.7 \psi1 \omega.5 \cal1 \end. }

\def\bar{\widebar}

\def\Abar{\,{\bar{\,\!A\,}}\,}




\phantom{}


\PrintVersionNumber



\sesquispace
\centerline{{\titlefont To What Type of Logic Does the ``Tetralemma'' Belong?}}



\bigskip


\singlespace			        

\author{Rafael D. Sorkin}
\address
 {Perimeter Institute, 31 Caroline Street North, Waterloo ON, N2L 2Y5 Canada}
\furtheraddress
 {Department of Physics, Syracuse University, Syracuse, NY 13244-1130, U.S.A.}
\email{rsorkin@pitp.ca}


\AbstractBegins                              
Although the so called tetralemma might seem to be incompatible with any
recognized scheme of logical inference, its four alternatives arise
naturally within the anhomomorphic logics which have been proposed in
order to accommodate certain features of microscopic (i.e. quantum)
physics.
This suggests that the possibility of similar, ``non-classical'' logics
might have been recognized in India at the time when Buddhism arose.
\AbstractEnds



\sesquispace
\vskip -30pt

\section{} 
Considered from the standpoint of classical logic, 
the fourfold structure of the so called tetralemma (\catuskoti)
appears to be irrational, 
and modern commentators have struggled to
explain its peculiar combination of alternatives, 
which at first sight appear to take 
the form, ``$A$'', ``not-$A$'', ``$A$ and not-$A$'', ``neither $A$ nor not-$A$''. 
(See for example [1] [2] [3] [4].)
Such a combination does not accord with contemporary logic,
and some authors have therefore concluded that the tetralemmas hail
from an unknown logical system which rejects one or more cherished
principles, like the ``Law of the excluded middle'' or 
the ``Law of contradiction''.\footnote{$^\star$}
{See Appendix II.}
Others commentators have abandoned any attempt to make logical sense of
the tetralemma, resorting rather to the hypothesis that its purpose is
merely to illustrate that the propositions that occur in it are
meaningless or ill-defined, like  ``The unicorn has
black eyes'', or (synesthesia aside) ``The number seven is yellow''.

According to Ruegg [3] the \catuskoti form occurs frequently in
the early Buddhist literature, and he provides (in translation) several
examples, one being the question whether the world is
\QuotationBegins
 finite, infinite, both finite and infinite, or neither finite nor infinite,
\QuotationEnds
another being the question whether
\QuotationBegins
 the world (of living beings) is eternal, not
  eternal, both eternal and not eternal, neither eternal nor not eternal.
\QuotationEnds
Similarly, Shcherbatskoi [1] writes that, according to
tradition, the founder of Buddhism refused to answer 
\QuotationBegins
four questions regarding the beginning of the world, viz., there is a
beginning, there is not, both, or neither.
\QuotationEnds
Almost the same
fourfold structure appears in another quotation from [3] (page 3):
\QuotationBegins
Entities of any kind are not ever found anywhere produced from
themselves, from another, from both [themselves and another], and also
from no cause.
\QuotationEnds
In all of these examples, the four offerings are evidently intended as
an exhaustive set of alternatives.  And with the possible exception of
the final tetralemma concerning causation, they all derive
from a single ``event'' or ``proposition'', $A$.
In the second example, for instance, $A$ is the event
of the world lasting forever, and if we abbreviate 
its negation
as $\Abar$, then it seems as if we could express the four
proffered alternatives as four propositions representable as:
$A\,$; \ 
$\Abar$; \ 
$A$ and $\Abar$; \ 
neither $A$ nor $\Abar$.\ 
The difficulty of course is that the final two alternatives are
self-contradictory when interpreted this way, leading one to wonder why
they had to be included at all.
(The tetralemma concerning causation might fall into a different
category than the other three, depending on whether or not one assumes
that every ``entity'' needs a cause or can have multiple causes.
Nevertheless, it still clashes with classical logic in its denial that
any one of the alternatives obtains.  This case is discussed further in
Appendix I.)

It seems clear that either the tetralemmas are unadorned nonsense or
something other than classical logic is involved.  In particular, one
cannot explain away the evident contradictions merely by hypothesizing
that questions about the beginning of the world or the infinity of
existence are meant to be like questions about unicorns or about the
color of the number seven.
Even if that
were the case, why would it not have sufficed to state the two
alternatives, $A$ and $\Abar$, after which one could either refuse to
choose between them or deny both of them?  Why did they need to be
supplemented with two further alternatives, apparently trivial or
lacking meaning altogether?

With respect to classical logic, this would have added nothing.
But could it be that the questioners were aware
of broader, more general logics, with respect to which the first two
alternatives failed to cover all the relevant possibilities, and with
respect to which the two additional alternatives were not in fact nonsense?
If so, what might these other logics have been?

A possible answer comes from quantum mechanics, where the paradoxes that
arise in the attempt to comprehend subatomic reality have given rise to
non-classical logical systems that aim to reflect more faithfully the nature of
the quantum world.
In the early proposals of this sort, known collectively as ``quantum
logic'', the laws for combining propositions via the connectives `and',
`or', and `not', are modified in such a way that the distributive law no
longer holds [5] [6] [7].
More recently, a different type of logical framework has been
put forward  
which modifies, not the connectives, but the rules of
inference [8] [9] [10].
It is these {\it\/anhomomorphic\/} logics, I would suggest, that hold the key
to understanding the \catuskoti form.

In order to appreciate 
why this might be,
one needs
to distinguish carefully between what have been called ``asserted'' and
``unasserted'' propositions.  That is, one needs to avoid confusing a
proposition-in-itself with the affirmation or denial of that
proposition.\footnote{$^\dagger$}
{In [11] Russel discusses this distinction at some length,
 contrasting his own attitude toward it with that of Frege.}
Unaccompanied by either affirmation or denial, a proposition 
functions rather like a question or a predicate.  
It only indicates a possible event or state of affairs, 
without taking a stand on whether 
that event or state of affairs 
actually happens or obtains.
For example, let $A$ be the proposition ``It rained all day yesterday''.
In itself, 
$A$ tells us nothing about yesterday's weather;
it only raises an implicit question.  
But if we then either deny or affirm $A$, we answer the implicit
question, either negatively or positively, as the case may be.

In ordinary speech we rarely if ever make this distinction explicitly,
though something rather like it is implicit in the grammar of
conditional and counterfactual constructions.  Perhaps for this reason,
formal logic also seems to lack an agreed upon symbolism to portray the
distinction.
%
To remedy this lack, we can (following [8]) introduce a
function $\phi$ that expresses affirmation or denial explicitly.  Given a
proposition\footnote{$^\flat$}
{{From the standpoint of physics, the word ``proposition'' could be
 misleading, because it could suggest that a formula like $\phi(A)=1$
 refers more to a sentence in a language than to the real world.  For
 this reason the word ``event'' (as in the world having a beginning) or
 ``state of affairs'' (as in the world being spatially finite) would be
 more apt than ``proposition''.  Nevertheless, I have retained the
 latter word in deference to the usage that is traditional within formal
 logic.} \ 
 In the nomenclature of [8] the function $\phi$ is
 called a {\it coevent}.  Using the ``turnstile'' notation, one could
 rewrite $\phi(A)=1$ as $\vdash{A}$, but there does not seem to be an
 analogous notation for $\phi(A)=0$.}
$A$, we can write $\phi(A)=1$ in order to {\it\/affirm\/} $A$, and
$\phi(A)=0$ in order to {\it\/deny\/} it.

The distinction between a proposition
per se and its affirmation or denial  is closely related to the
distinction between what, using a different language, might have been 
called ``positive'' and ``negative'' propositions.  
Consider the (unasserted) propositions,
$A$ = ``the electron is here''
and
$\Abar$ = ``the electron is elsewhere''.
A statement like ``I see the electron here'' is in some sense positive,
corresponding roughly to $\phi(A)=1$ (i.e. to the affirmation of $A$).
Similarly, ``I see the electron there'' 
is also positive, and
corresponds roughly to $\phi(\Abar)=1$.
In contrast, a statement like
``I don't see the electron here'' is in
some sense negative, 
corresponding roughly to $\phi(A)=0$ 
(i.e. to the denial of $A$),
while 
 ``I don't see the electron there'' 
corresponds roughly to $\phi(\Abar)=0$.
Viewed in this manner, the affirmation of $A$ is decoupled from the denial of
$\Abar$, and four distinct alternatives arise,
corresponding to those of the tetralemma.

Physically, the motivation for this type of de-coupling stems from the
phenomenon of {\it interference}, in which, for example, an electron can
act as if it were ``in two places at once'' or ``in neither place
separately but both together'' [12].  Anhomomorphic logics
accommodate this type of behavior by admitting, for instance, the
possibility that both of the
propositions, ``the electron is here'' and ``the
electron is elsewhere'' might be false.  In other words, one can have
$\phi(A)=\phi$($\Abar$)=0 in such logics, just as one can have
$\phi(A)=\phi$($\Abar$)=1.\footnote{$^\star$}
{Not every scheme of anhomomorphic logic recognizes all of these
 possibilities, however.  For example in the {\it multiplicative scheme}
 described in [8] the combination
 $\phi(A)=\phi$($\Abar$)=1 cannot occur.}
\ Thus, there is no necessary correlation in anhomomorphic
logic between $\phi(A)$ and $\phi(\hbox{not-$A$})$.

Consider now the sentence ``The electron is not elsewhere'' (or if you
wish, ``The world is not infinite'').  With the possibility of an
anhomomorphic logic in mind, we can see that its meaning is ambiguous.
Perhaps it is simply formulating an unasserted proposition-in-itself, namely the
above proposition $A$ = ``the electron is here''.  
On the other hand, perhaps it is trying to {\it deny} the complementary
proposition $\Abar$, 
i.e. to express that $\phi(\Abar)=0$.  
Or perhaps it actually means
to {\it affirm} the proposition $A$, 
in which case one should render it as $\phi(A)=1$.  
Accordingly, three different readings of our sentence are possible: 
$(i)$ $A$, $(ii)$ $\phi(\Abar)=0$, $(iii)$ $\phi(A)=1$.
In classical logic, there is no good reason to distinguish these
meanings, but anhomomorphically it is crucial to do so.  I would
conjecture that most of the confusion over the tetralemma can be traced
to a wrong resolution of this unrecognized ambiguity.

Take for definiteness the \catuskoti about whether the world is finite,
letting $A$ be the proposition ``the world is finite'', and $\Abar$ the
proposition ``the world is infinite'', which is the logical negation of
$A$.\footnote{$^\dagger$}
{It is is not clear that an analysis in terms of complementary
 propositions $A$ and $\Abar$ applies to the tetralemma about causation.
 See the appendices on this and related issues.}
The question is how to interpret the alternatives that
constitute the tetralemma. 

If we interpret them as propositions-per-se, we will immediately
land back in the usual difficulties, since even within anhomomorphic
logic, the propositions,\ 
``both $A$ and $\Abar\,$'' 
\ and \ 
``neither $A$ nor $\Abar\,$'' \  
are still nonsensical.  
(More precisely,  both of them reduce to the ``zero proposition'', and including
them as distinct alternatives adds nothing, since they could never be
chosen as the answer.)
Moreover, we can gain no comfort from any ambiguity in the symbolic form
of these propositions, since for example, the
propositions,
\hbox{$\Abar\wedge\widebar{\Abar}$} 
and
$\overline{A\vee\Abar}$ 
are strictly equivalent
within anhomomorphic logic (and both equal to zero).\footnote{$^\flat$}
{Here, I've used for brevity the common symbolism, $\vee$=`or', $\wedge$=`and'.}
\ (For this reason, the proposal in [4] to distinguish between
these two forms of the fourth alternative would get nowhere within
anhomomorphic logic.)
%

However if --- in resolving the ambiguity pointed out above ---
we read the four alternatives differently, namely as the four possible
combinations of affirmation and denial of the two complementary
propositions $A$ (``the world is finite'') and $\Abar$ (``the world is infinite''), 
then the tetralemma makes perfect sense.  
In symbols, its four alternatives become then:

(1) $\phi(A)=1$ and $\phi(\Abar)=0$ 

(2) $\phi(A)=0$ and $\phi(\Abar)=1$

(3) $\phi(A)=1$ and $\phi(\Abar)=1$

(4) $\phi(A)=0$ and $\phi(\Abar)=0$

\noindent
On this exegesis, the tetralemma format is simply the one that you adopt
if you are seeking maximum generality, once you have admitted that the
denial of a proposition $A$ need not entail the affirmation of its logical
complement $\Abar$, and vice versa.
To the extent that Indian thinkers in the time of Gotama were aware of
this possibility, they would naturally have phrased their questions in
``tetralemmatic'' form.  (And of course, tradition would ensure that
the form would persist, even if its original raison d'etre had been lost.)

But do we possess independent evidence that they were in fact aware of such
more general forms of logic?
Perhaps attention to the precise wording of the various \catuskoti in
their original languages would shed light on this question, but we
also have at least one indication which 
appears to be relatively
independent of
that kind of delicate textual analysis.
In reference [3], Ruegg writes that Indian grammarians recognized
a type of ``absolute'' or ``pure'' negation that did not entail
``affirmation of the contrary''.  
Could this type of negation be what we have been expressing as $\phi(A)=0$,
and to which one should oppose 
the ``propositional negation'' that exchanges $A$ with $\Abar$?
And if so,
could the introduction of two types of negation
have been their way to decouple the denial of $A$
from the affirmation of $\Abar$, and vice versa?
\footnote{$^\star$}
{Westerhoff [13] also emphasizes the relevance of the two types of
 negation, {\it\/prasajya\/} and {\it\/\paryudasa\/}.  In my
 interpretation, the former, or ``absolute'' negation corresponds to the
 denial of $A$, $\phi(A)=0$, while the latter corresponds to the
 ``propositional negation'', $A \longleftrightarrow \Abar$.}

It would be very interesting to know what led people to
recognize --- if they did recognize --- the possibility of an
anhomomorphic logic over two millenia ago.
They cannot have had access to the kind of technology that has led in
modern times to quantum physics.  Are there then other experiences that one
could point to which were in fact available to them and to which
anhomomorphic inference is more suited than homomorphic inference?  If
so, we might gain a better intuition for the microworld by ourselves
paying more attention to those experiences.

\section{Appendix I. A second-order application of anhomomorphic logic?}
The analysis given in the body of the paper applies straightforwardly to
all but the last example cited at the outset (about causation).  If
you accept it, it explains why these questions were posed in fourfold
form.  But if that were the whole story, then the answer to any given
tetralemma should be a selection of one of the four mutually exclusive
and exhaustive alternatives,
(1)$-$(4), as the correct one.  
According to [13], the answers have indeed been like that in many
cases.  But in other cases, all four have been accepted or all four
rejected, notably in the work of \Nagarjuna, who (unlike Gotama, who
refused to answer at all) in his writings rejected all four of the
alternatives [3].


Perhaps, as is commonly suggested, \Nagarjuna was simply trying to 
express a mystical rejection of analytical thought itself.  However,
it seems worth pointing out that anhomomorphic logic opens up another
interpretation, perhaps consistent with the mystical one, but not really
requiring it.  Namely one can imagine that \Nagarjuna's blanket denial
represents a kind of ``second order'' application of anhomomorphic
logic, one that reasons anhomomorphically, not just about ``reality as
such'', but also about the logical processes employed to grasp that
reality.  

When one is dealing with two propositions $A$ and $B$ (whether
or not they are mutual negations so that $B=\Abar$), there will be
four possible combinations of affirmation and denial of each, as
explained above, i.e. four possible combinations of $\phi(A)$ and
$\phi(B)$.  
When one is dealing with not two but four propositions,  
there will
be 16 possible combinations, including the one that denies all four.
\Nagarjuna  could be 
opting for this last combination
in respect of the four
``second-order propositions'' listed above as (1)$-$(4).  It is natural to
regard them as of second order since they are 
in some sense ``propositions about propositions'', speaking not just
about ``first order events'', but about the values of $\phi$ on these
``events''.  
Thus, for example, one might try to symbolize the denial of alternative (3)
by writing  
$\phi(\,\phi(A)=\phi(\Abar)=1\,) \,=\, 0$.

Whether the tetralemma on causation also involves such a ``denial at
second order'' seems uncertain.  On one hand, it can be fit into the
same mold as the others, if one treats ``caused by self'' and ``caused
by another'' as mutual negations, and if one interprets its negative
phrasing as the denial of all four alternatives.  In that case, it would
have the same character as the other examples, and Nargarjuna's blanket
denial could again be interpreted either in the mystical or the ``second
order'' manner.

However, such an exegesis would only be correct if one assumed that
every event (or ``entity'' in the quoted translation) must have a unique
cause.  Without that assumption, the 
propositions $A$ = ``caused by self'' and $B$ = ``caused by another''
are not logical complements of each other, because the event in question
might have multiple causes or no cause at all.
In that case, it is not necessary (although it would still be possible)
to interpret \Nagarjuna's blanket denial in either of the above manners,
because the naive symbolic-logic reading of the four alternatives is no
longer problematic.  That is, it could be that the alternatives are
merely the perfectly respectable propositions: 
$A$; $B$; $A$ and $B$; neither $A$ nor $B$.  Then one would have in this
tetralemma a set of alternatives which happen to be four in number, not
because they are formed from a single proposition $A$ in the manner of
the other tetralemmas, but because they represent all possible subsets
of the two-element 
set \lbrace``caused by self'', ``caused by another''\rbrace.  In that
case, denying all the alternatives would only be a ``first order''
(though still anhomomorphic) activity, symbolizable as
$\phi(A)=\phi(B)=\phi(A\wedge B)=\phi(\hbox{not}(A\vee B))=0$.

\section{Appendix II. Contradiction and excluded middle}
It is often felt that the tetralemma's phrasing conflicts with the
so-called laws of contradiction and of the excluded middle.  Is
such a conflict still present on the anhomomorphic interpretation?  If
not, this would furnish further evidence for the interpretation,
inasmuch as --- far from denying these two laws of logic --- prominent
schools of Buddhist logic seem to have embraced them, in particular the
Madhyamaka school [3].
  
What these laws mean in modern terms seems to be not exactly settled,
but to the extent that they govern the formation of
compound propositions from simpler ones, they hold automatically in
anhomomorphic logic.  For example not-not-$A$ = $A$ holds for any
proposition $A$.

To the extent, however, that they are understood as laws of inference,
the possibility of conflict does arise.  For present purposes, let us
take ``excluded middle'' to mean that 0 (false) and 1 (true) are the
only values a proposition can assume, and let us take
``non-contradiction'' to mean that no proposition can take both of these
values at once.  Thus interpreted, both laws are honored in
anhomomorphic logic, at least as it has developed so far.  The former is
obeyed because $\phi(A)=0$ and $\phi(A)=1$ are the only two values
provided for, the latter simply because $\phi$ is assumed to be a
function, as opposed to what used to be called a ``multiply valued
function''.



\bigskip
\noindent
Research at Perimeter Institute for Theoretical Physics is supported in
part by the Government of Canada through NSERC and by the Province of
Ontario through MRI.

\ReferencesBegin                             

\ref [1] Th. Stcherbatsky, {\it Buddhist Logic}, volume I (Dover 1962), page 477

\ref [2] A.K. Warder, ``Dharmas and data'',
   \journaldata{Journal of Indian Philosophy}{1(3)}{272-295}{1971}

\ref [3] D. Seyfort Ruegg, ``The uses of the four positions of the \catuskoti
   and the problem of the description of reality in Mah{\=a}y{\=a}na Buddhism'',
   \journaldata{Journal of Indian Philosophy}{5}{1-71}{1977}

\ref [4] Sitansu S. Chakravarti, ``The M{\=a}dhyamika \catuskoti 
   or Tetralemma'',
   \journaldata{Journal of Indian Philosophy}{8}{303-306}{1980}

\ref [5] G. Birkhoff and J. von Neumann, ``The Logic of Quantum Mechanics'', 
  \journaldata{Annals of Mathematics}{37}{823-843}{1936}

\ref [6]  H. Putnam, ``Is Logic Empirical?'', in 
  {\it Boston Studies in the Philosophy of Science vol. 5},
  M. Wartofsky and R. Cohen, eds. 
  (Humanities Press, New York, 1969), pp. 216-241

\ref [7] 
David Finkelstein, 
``Matter, Space and Logic'', 
in 
{\it Boston Studies in the Philosophy of Science vol. 5},
  M. Wartofsky and R. Cohen, eds. 
  (Humanities Press, New York, 1969), pp. 199-215

\ref [8] Rafael D. Sorkin, ``An exercise in `anhomomorphic logic'~'',
 \journaldata{Journal of Physics: Conference Series (JPCS)}{67}{012018}{2007},
 a special volume edited by L. Diosi, H-T Elze, and G. Vitiello, and
 devoted to the Proceedings of the DICE-2006 meeting,
  held September 2006, in Piombino, Italia.
  \eprint{arxiv quant-ph/0703276} , \lbr
  \eprint{http://www.pitp.ca/personal/rsorkin/some.papers/}

\ref [9] Yousef Ghazi-Tabatabai, ``Quantum Measure Theory: A New Interpretation'',
 Doctoral dissertation (Imperial College, London, January 2009)
 \lbr
 \arxiv{0906.0294}

\ref [10] Rafael D.~Sorkin, ``Does a quantum particle know its own energy?''
  \journaldata{Journal of Physics: Conf. Ser.}{442}{012014}{2013}
  \eprint{http://arxiv.org/abs/1304.7550}
  \lbr
  \eprint{http://www.pitp.ca/personal/rsorkin/some.papers/147.know.own.energy.pdf}

\ref [11] Bertrand Russel, {\it The Principles of Mathematics, second edition} 
(New York, W.W. Norton)

\ref [12] The experimental grounds for such descriptions
can be found in any number of textbooks.  For a particularly vivid
account see Section 1-5 in: \lbr
Richard P.~Feynman, Robert B.~Leighton and Matthew~Sands,
  {\it The Feynman Lectures on Physics, vol. III: Quantum Mechanics}
  (Addison--Wesley, 1965)

\ref [13] Jan Westerhoff, ``\Nagarjuna's \Catuskoti'', \journaldata{J. Indian Philosophy}{34}{367-395}{2006}

\end                                         


(prog1 'now-outlining
  (Outline* 
     "\f"               1
      "
      "
      "
      "\\Abstrac"       1
      "\\section"       1
      "\\subsectio"     3
      "\\appendi"       2
      "\\Referen"       1
      "\\ref "          3
      "\\end